\documentclass[12pt]{article}
\usepackage{amsmath,amssymb,theorem}
\usepackage{graphicx}
\usepackage{subfigure}
\usepackage{psfrag}
\usepackage{color}
\usepackage{enumerate}
\usepackage{epstopdf}
\usepackage[T1]{fontenc}
\usepackage[colorlinks,linkcolor=blue]{hyperref}
\newtheorem{thm}{Theorem}[section]
\newtheorem{conj}[thm]{Conjecture}
\newtheorem{cor}[thm]{Corollary}
\newtheorem{lem}[thm]{Lemma}
\theorembodyfont{\rmfamily}
\def\pf{\bigskip\noindent {\bf Proof.}~~}

\def\dfn#1{{\sl #1}}

\def\less{\backslash}

\def\pf{\bigskip\noindent {\bf{Proof.}}~~}

\def\fattextindent#1{\indent\indent\llap{#1\enspace}\ignorespaces}

\def\myitemitem{\par\hangindent\parindent\fattextindent}

\title{Star 5-edge-colorings of subcubic multigraphs}
\author{Hui Lei$^1$, Yongtang Shi$^1$, Zi-Xia Song$^2$\thanks{Corresponding author.  \newline Email addresses:  leihui0711@163.com (H. Lei), shi@nankai.edu.cn (Y. Shi), Zixia.Song@ucf.edu (Z-X. Song), wangtao@henu.edu.cn (T. Wang)}, Tao Wang$^3$ \\
{ $^1$ Center for Combinatorics and LPMC}\\
{ Nankai University, Tianjin 300071,  China}\\
{ $^2$ Department  of Mathematics}\\
{ University of Central Florida, Orlando, FL32816, USA}\\
{ $^3$ Institute of Applied Mathematics}\\
{ Henan University, Kaifeng, 475004, P.R. China}\\
}

\begin{document}
\maketitle
\begin{abstract}
The \dfn{star chromatic index}  of a multigraph $G$, denoted
$\chi'_{s}(G)$,  is the minimum number of colors needed to properly
color the edges of $G$ such that no path or cycle of length four is
bi-colored.  A multigraph $G$ is \dfn{star $k$-edge-colorable} if
$\chi'_{s}(G)\le k$.   Dvo\v{r}\'ak, Mohar and \v{S}\'amal [Star
chromatic index,  J. Graph Theory  72 (2013), 313--326]
proved  that every  subcubic multigraph is star $7$-edge-colorable, and  conjectured  that every  subcubic multigraph
should be star $6$-edge-colorable.       Kerdjoudj, Kostochka  and Raspaud  considered the list version of this problem for simple graphs and proved that every subcubic graph with maximum average degree  less than $7/3$ is  star list-$5$-edge-colorable. It is known that  a graph with maximum average degree $14/5$  is not necessarily  star $5$-edge-colorable. In this paper, we  prove that every subcubic multigraph  with maximum average degree  less than $12/5$ is star $5$-edge-colorable.  \\

\end{abstract}

\section{Introduction}

\baselineskip 17pt

All multigraphs  in this paper are finite and loopless; and all graphs  are finite and without loops or multiple edges.
Given a multigraph $G$, let $c: E(G)\rightarrow [k]$ be a proper
edge-coloring of $G$, where $k\ge1$ is an integer and $[k]:=\{1,2, \dots, k\}$. We say that
$c$ is a  \dfn{star $k$-edge-coloring} of $G$ if no path or cycle of
length four in $G$ is bi-colored under the coloring $c$; and  $G$ is
\dfn{star $k$-edge-colorable} if $G$ admits a star
$k$-edge-coloring. The \dfn{star chromatic index}  of $G$, denoted $\chi'_{s}(G)$,  is the smallest integer $k$ such that $G$ is
star $k$-edge-colorable.  
As pointed out in  \cite{DMS2013}, the definition of
star edge-coloring of a graph $G$ is equivalent to the star
vertex-coloring of its line graph $L(G)$.  Star edge-coloring of a
graph was initiated  by Liu and Deng \cite{DL2008}, motivated by the
vertex version (see \cite{ACKKR2004, BCMRW2009, CRW2013, KKT2009,
NM2003}).   Given a multigraph $G$, we use  $|G|$ to denote the number of vertices,  $e(G)$ the number of edges, $\delta(G)$ the minimum degree, and $\Delta(G)$ the maximum degree of $G$, respectively.  
We use $K_n$ and $P_n$ to denote the complete graph and the path on $n$ vertices, respectively.  A multigraph $G$ is \dfn{subcubic} if all its vertices have degree less than or equal to
three. The \dfn{maximum average degree} of  a multigraph $G$, denoted  $\text{mad}(G)$, is defined as the maximum  of  $2
e(H)/|H|$ taken over all the subgraphs $H$ of $G$.  The following upper bound is a result of  Liu and Deng \cite{DL2008}.  \medskip

\begin{thm}[\cite{DL2008}]
For any graph  $G$ with  $\Delta(G)\geq7$,
$\chi'_{s}(G)\leq \lceil16(\Delta(G)-1)^\frac3 2\rceil.$
\end{thm}

Theorem~\ref{Kn}  below is a result of Dvo\v{r}\'ak, Mohar and \v{S}\'amal~\cite{DMS2013},   which gives an   upper and a lower bounds for complete graphs.

\begin{thm} [\cite{DMS2013}]\label{Kn}  The star chromatic index of   the complete graph $K_n$ satisfies

$$2n(1+o(1))\leq \chi'_{s}(K_n)\leq n\, \frac{2^{2\sqrt{2}(1+o(1))\sqrt{\log n}}}{(\log n)^{1/4}}.$$
In particular, for every $\epsilon>0$, there exists a constant $c$ such that  $\chi'_{s}(K_n)\le cn^{1+\epsilon}$ for every integer $n\ge1$.
\end{thm}

The true order of magnitude of  $\chi'_{s}(K_n)$ is still unknown. From Theorem~\ref{Kn},  an upper bound in terms of the maximum degree for general  graphs is also derived in~\cite{DMS2013}, i.e.,
$\chi'_{s}(G)\leq \Delta\cdot 2^{O(1)\sqrt{\log \Delta}}$ for  any  graph $G$ with maximum degree $\Delta$.  In the same paper, Dvo\v{r}\'ak, Mohar and \v{S}\'amal~\cite{DMS2013} also considered the star chromatic index of subcubic multigraphs. To state their result, we need to introduce one notation. A graph $G$ \dfn{covers} a graph $H$ if there is a mapping $f: V(G)\rightarrow V(H)$ such that for any $uv\in E(G)$, $f(u)f(v)\in E(H)$, and for any $u\in V(G)$, $f$ is a bijection between $N_G(u)$ and $N_{H}(f(u))$.  They proved the following.

\begin{thm} [\cite{DMS2013}]\label{s=7} Let $G$ be a multigraph.
\begin{enumerate}[(a)]
\item  If $G$ is  subcubic,  then $\chi'_s(G)\le7$.\vspace{-8pt}

\item   If $G$ is  cubic and has no multiple edges, then $\chi'_s(G)\ge4$ and the equality holds if and only if $G$ covers the graph of $3$-cube.
\end{enumerate}
\end{thm}

As observed in~\cite{DMS2013},  $K_{3,3}$ is  not star $5$-edge-colorable but star $6$-edge-colorable. No subcubic multigraphs with star chromatic index seven  are known. Dvo\v{r}\'ak, Mohar and \v{S}\'amal~\cite{DMS2013}   proposed the following conjecture.

\begin{conj} [\cite{DMS2013}]\label{cubic}
Let $G$ be a subcubic multigraph. Then $\chi'_s(G)\leq 6$.
\end{conj}

  It was  shown in~\cite{BLM2016} that every subcubic outerplanar graph is star $5$-edge-colorable.
 Lei, Shi and Song~\cite{LSS2017} recently proved that every subcubic multigraph $G$ with $\text{mad}(G)<24/11$ is star $5$-edge-colorable, and  every subcubic multigraph $G$ with $\text{mad}(G)<5/2$ is star $6$-edge-colorable. Kerdjoudj,  Kostochka and Raspaud~\cite{KKP2017} considered the list version of star edge-colorings of simple graphs. They proved that every subcubic graph is   star list-$8$-edge-colorable, and further proved the following stronger results.
\medskip
\begin{thm} [\cite{KKP2017}]\label{KKP} Let $G$ be a subcubic graph.
\begin{enumerate}[(a)]
\item If $\text{mad}(G)<7/3$, then $G$ is  star list-$5$-edge-colorable.\vspace{-10pt}
\item If $\text{mad}(G)<5/2$, then $G$ is  star list-$6$-edge-colorable.
\end{enumerate}
\end{thm}

 \begin{figure}[htbp]
\begin{center}
\includegraphics[scale=0.4]{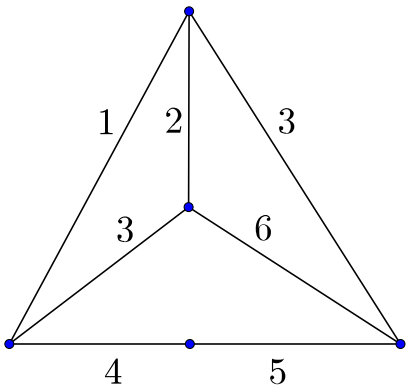}
\caption{ A graph  with maximum average degree $14/5$ and star chromatic index  $6$.}\label{prism}
\end{center}
\end{figure}
As mentioned above,   $K_{3,3}$ has star chromatic index $6$, and is  bipartite and non-planar.  The graph, depicted in Figure~\ref{prism},  has star chromatic index $6$, and is planar and  non-bipartite.  We see that not every bipartite, subcubic graph is star $5$-edge-colorable; and not every planar, subcubic graph is star $5$-edge-colorable. It remains unknown whether every bipartite, planar subcubic multigraph is star $5$-edge-colorable.  
In this paper, we improve Theorem~\ref{KKP}(a) by showing the following main result.

\begin{thm}\label{mainthm}  Let $G$ be a subcubic multigraph with $mad(G)<12/5$. Then $\chi'_{s}(G)\leq 5$.
\end{thm}

  We don't know if the bound $12/5$  in Theorem~\ref{mainthm} is best possible.  The graph depicted  in Figure~\ref{prism}  has  maximum average degree $14/5$  but is not star $5$-edge-colorable.  \medskip

The  \dfn{girth} of a graph $G$ is the length of a shortest cycle in $G$.  It was observed in~\cite{girth} that  every planar graph with girth $g$ satisfies $\text{mad}(G)< \frac {2g}{g-2}$.   This, together with Theorem~\ref{mainthm},  implies the following.

\begin{cor}
Let $G$ be a planar subcubic graph with girth $g$. If $g\geq12$, then $\chi'_{s}(G)\leq 5$.\medskip
\end{cor}

 We need to introduce more notation. Given a multigraph $G$,
 a vertex of degree $k$  in $G$ is a \dfn{k-vertex}, and a \dfn{k-neighbor} of a vertex $v$  in $G$ is a $k$-vertex adjacent to $v$ in $G$.  A \dfn{$3_k$-vertex} in $G$
is a $3$-vertex  incident to exactly $k$ edges $e$ in $G$ such that the other end-vertex  of $e$ is a $2$-vertex. For any proper edge-coloring  $c$ of a multigraph $G$ and for any $u\in V(G)$, let $c(u)$ denote the set of all colors such that each is used to color an edge incident with $u$ under the coloring $c$.
For any two sets  $A, B$,   let  $A\less B := A-B$.   If  $B=\{b\}$, we simply write  $A\less b$ instead of $A\less B$. \medskip




\section{Properties of star $5$-critical subcubic multigraphs }\label{prop}

A  multigraph $G$ is  \dfn{star $5$-critical} if $\chi'_s(G)>5$ and  $\chi'_s(G-v)\le 5$ for any $v\in V(G)$. In this section,  we establish some structure results on star $5$-critical subcubic multigraphs.   Clearly, every star $5$-critical multigraph  must be connected.
\medskip

Throughout the remainder  of this section, let $G$ be a star $5$-critical subcubic multigraph, 
and let $N(v)$ and $d(v)$ denote the neighborhood and  degree of a vertex $v$ in $G$, respectively.  Since every multigraph with maximum degree at most  two or  number of vertices at most four  is star $5$-edge-colorable, we see that  $\Delta(G)=3$  and $|G|\ge5$.  As observed in \cite{LSS2017}, any $2$-vertex in $G$ must have  two distinct neighbors.  The following Lemma~\ref{deg=1} and Lemma~\ref{deg=2} are proved in \cite{LSS2017} and will be used in this paper.   

\begin{lem}[\cite{LSS2017}]\label{deg=1} For any $1$-vertex $x$ in $G$,  let  $N(x)=\{y\}$. The following are true.
\begin{enumerate}[(a)]
\item   $|N(y)|=3$. \vspace{-10pt}
\item $N(y)$ is an independent set in $G$,  $d(y_1)=3$ and $d(y_2)\ge 2$, where $N(y)=\{x, y_1,y_2\}$ with   $d(y_1)\ge d(y_2)$.\vspace{-10pt}
\item If $d(y_2)=2$, then for any $i\in\{1,2\}$ and any $v\in N_G(y_i)\less y$,   $|N(v)|\ge2$, $|N(y_1)|=3$, $|N(y_2)|=2$,  and $N[y_1]\cap N[y_2]=\{y\}$.\vspace{-10pt}
\item  If $d(y_2)=2$, then  $d(w_1)=3$, where $w_1$ is the other neighbor of $y_2$ in $G$.\vspace{-10pt}
\item If  $d(y_2)=3$, then either $d(v)\ge2$ for any $v\in N(y_1)$ or $d(v)\ge2$ for any $v \in N(y_2)$.
\end{enumerate}
\end{lem}

 \begin{lem}[\cite{LSS2017}]\label{deg=2}
For any $2$-vertex $x$ in $G$,  let   $N(x)=\{z, w\}$ with $|N(z)|\le |N(w)|$. The following are true.
\begin{enumerate}[(a)]
\item   If $zw\in E(G)$, then   $|N(z)|=|N(w)|=3$ and $d(v)\ge2$ for any $v\in N(z)\cup N(w)$.\vspace{-10pt}
\item  If  $zw\notin E(G)$, then  $|N(w)|=3$ or  $ |N(w)|=|N(z)|=2$,  and $d(w)=d(z)=3$.\vspace{-10pt}
\item   If $d(z)=2$ and  $z^*w\in E(G)$,  then   $|N(z^*)|=|N(w)|=3$, and $d(u)=3$ for any $u\in (N[w] \cup N[z^*])\less \{x,z\}$, where $z^*$ is the other neighbor of $z$ in $G$.\vspace{-10pt}
\item  If 
$d(z)=2$, then  $|N(z^*)|=|N(w)|=3$, and $|N(v)|\geq2$ for any $v\in N(w)\cup N(z^*)$, where $N(z)=\{x, z^*\}$.
\end{enumerate}
\end{lem}
\medskip

Let $H$ be the graph obtained from $G$ by deleting all $1$-vertices.     By Lemma~\ref{deg=1}(a,b), $H$ is connected and  $\delta(H)\geq2$.  Throughout the remaining of the proof,  a $2$-vertex   in $H$ is \dfn{bad} if  it has a  $2$-neighbor  in $H$, and a $2$-vertex   in $H$ is \dfn{good} if  it is not bad.  For any $2$-vertex $r$ in $H$, we use $r'$ to denote the unique $1$-neighbor of $r$ in $G$ if $d_G(r)=3$. By Lemma~\ref{deg=1}(a) and the fact that  any $2$-vertex in $G$ has two distinct neighbors in $G$, we obtain the following two lemmas.

 \begin{lem}\label{2nbr}
 For any $2$-vertex $x$ in $H$,  $|N_{H}(x)|=2$.
 \end{lem}

 \begin{lem}\label{3nbr}
 For any $3_k$-vertex   $x$ in $H$ with $k\ge2$,  $|N_{H}(x)|=3$.
 \end{lem}
\medskip

Proofs of Lemma~\ref{noC3} and Lemma~\ref{noC4} below can be obtained  from the proofs  of Claim 11  and   Lemma 12 in~\cite{KKP2017}, respectively. Since a star $5$-critical multigraph is not necessarily the edge minimal counterexample in the proof of  Theorem 4.1 in~\cite{KKP2017},   we include new proofs of  Lemma~\ref{noC3} and Lemma~\ref{noC4} here for completeness.

 \begin{lem}\label{noC3}
$H$ has no $3$-cycle such that  two of its  vertices  are bad.
\end{lem}

\pf Suppose that $H$ does contain a $3$-cycle with vertices $x, y, z$   such that both $y$ and $z$ are bad.  Then $x$ must be a $3$-vertex  in $G$ because $G$ is $5$-critical. Let $w$ be the third neighbor of $x$ in $G$. Since $G$ is $5$-critical,  let $c: E(G\less \{y,z\})\rightarrow [5]$ be any star $5$-edge-coloring of $G\less \{y,z\}$.  Let $\alpha$ and $\beta$ be two distinct numbers in $[5]\less c(w)$ and $\gamma\in [5]\less \{\alpha, \beta, c(xw)\}$. Now coloring the edges $xy, xz, yz$  by colors  $\alpha, \beta, \gamma $  in order,  and further  coloring all   the edges $yy', zz'$ by color $c(xw)$ if $y'$ or  $z'$ exists,  we obtain a star $5$-edge-coloring of $G$, a contradiction. 
 \hfill\vrule height3pt width6pt depth2pt\medskip

 \begin{lem}\label{noC4}   $H$ has no $4$-cycle  with vertices  $x,u,v,w$ in order  such that all of $u,v,w$ are  bad. Furthermore, if $H$ contains a path with vertices $x,u,v,w,y$ in order such that all of  $u, v,w$ are  bad, then both $x$ and $y$ are $3_1$-vertices in $H$.
\end{lem}

\pf 
Let $P$ be a path   in $H$ with  vertices $x,u,v,w,y$ in order such that all of  $u, v,w$ are  bad, where $x$ and $y$ may be the same.  Since all of  $u, v,w$ are bad, by the definition of $H$, $uw\notin E(G)$.
  By Lemma~\ref{deg=1}(b,c,e) applied to the vertex $v$,  $d_G(v)=2$.  By Lemma~\ref{deg=2}(b) applied to $v$,   $d_G(u)=d_G(w)=3$. Thus both $w'$ and  $u'$ exist.  Now by Lemma~\ref{deg=1}(c) applied to $u'$ and $w'$,  $d_{H}(x)=d_{H}(y)=3$, and $x\ne y$. This proves  that $H$ has no $4$-cycle  with vertices  $x,u,v,w$ in order  such that all of $u,v,w$ are  bad.  \medskip

We next show that both $x$ and $y$ are $3_1$-vertices in $H$. Suppose that one of $x$ and $y$, say $y$, is not a   $3_1$-vertex in $H$.  Then  $y$ is either a $3_2$- vertex or $3_3$-vertex   in $H$.  By Lemma~\ref{3nbr},  $|N_{H}(y)|=3$.  Let $N_{H}(y)=\{w,y_1,y_2\}$ with $d_{H}(y_1)=2$. Then  $y_1\neq u$, otherwise $H$ would have a $4$-cycle with vertices  $y,u,v,w$ in order  such that all of $u,v,w$ are  bad.  Note that $y_2$ and $x$ are not necessarily distinct.
By Lemma~\ref{2nbr}, let  $r$ be the other neighbor of $y_1$ in $H$. Since $G$ is $5$-critical,  let $c: E(G\less \{v,u',w'\})\rightarrow [5]$ be any star $5$-edge-coloring of $G\less \{v,u',w'\}$. We may assume that $c(wy)=3$, $c(yy_1)=1$ and $c(yy_2)=2$. We first color  $uv$ by a color $\alpha$ in $[5]\less (c(x)\cup\{3\})$ and $uu'$ by a color $\beta$ in $[5]\less (c(x)\cup\{\alpha\})$.
Then $3\in c(y_1)\cap c(y_2)$, otherwise, we may assume that $3\notin c(y_i)$ for some $i\in\{1,2\}$, now coloring $vw$ by a color $\gamma$ in $\{i,4,5\}\less\alpha$ and  $ww'$ by a color in $\{i,4,5\}\less\{\alpha,\gamma\}$
yields a star $5$-edge-coloring of $G$, a contradiction. It follows that  $4, 5\in c(y_1)\cup c(y_2)$, otherwise, say $\theta\in\{4,5\}$ is not in $c(y_1)\cup c(y_2)$, now recoloring $wy$ by color $\theta$, $uv$ by a color $\alpha'$ in $\{\alpha,\beta\}\less\theta$, $uu'$ by $\{\alpha,\beta\}\less\alpha'$, and then coloring $ww'$ by a color in $\{1,2\}\less\alpha'$ and $vw$ by a color in $\{3,9-\theta\}\less\alpha'$, we obtain a star $5$-edge-coloring of $G$, a contradiction.  Thus $c(y_1)=\{1,3,\theta\}$ and   $c(y_2)=\{2,3,9-\theta\}$, where $\theta\in\{4,5\}$.  If $c(y_1y_1')\neq3$ or $c(y_1r)=\theta$ and $1\notin c(r)$, then we obtain a star $5$-edge-coloring of $G$ by recoloring $wy$ by color $\theta$, $uv$ by a color $\alpha'$ in $\{\alpha,\beta\}\less\theta$,  $uu'$ by $\{\alpha,\beta\}\less\alpha'$, and then coloring $ww'$ by a color $\gamma$ in $\{2,3,9-\theta\}\less\alpha'$,  and $vw$ by a color in $\{2,3,9-\theta\}\less\{\alpha',\gamma\}$.
Therefore,  $c(y_1y_1')=3$ and $1\in c(r)$. Now recoloring $y_1y_1'$ by a color in $\{2,9-\theta\}\less c(r)$,  we obtain a star $5$-edge-coloring $c$ of $G\less \{v,u',w'\}$ satisfying $c(wy)=3$, $c(yy_1)=1$   and $c(yy_2)=2$ but $3\notin c(y_1)\cap c(y_2)$, a contradiction.  Consequently,  each of  $x$ and $y$ must be a  $3_1$-vertex in $H$. This completes the proof of Lemma \ref{noC4}.  \hfill\vrule height3pt width6pt depth2pt\medskip

 \begin{lem}\label{noBad}  For any $3_3$-vertex $u$ in $H$,  no vertex in $N_{H}(u)$ is  bad.
 \end{lem}

 \pf Let $N_{H}(u)=\{ x, y, z\}$ with  $d_{H}(x)=d_{H}(y)=d_{H}(z)=2 $. By Lemma~\ref{3nbr}, $u,x,y,z$ are all distinct.    By Lemma~\ref{2nbr}, let $x_1$,  $ y_1$ and $z_1$ be  the other neighbors of $x, y, z$ in $H$, respectively. Suppose that  some vertex, say $x$, in $N_{H}(u)$  is bad. Then $d_{H}(x_1)=2$.  By Lemma~\ref{2nbr}, let $w$ be the other neighbor  of $x_1$ in $H$.
 By Lemma~\ref{noC3} and Lemma~\ref{noC4}, $N_{H}(u)$ is an independent set and $x_1\notin \{y, z,y_1, z_1\}$.    Notice that  $y_1$,  $z_1$  and $w$  are not necessarily distinct.   Let $A:=\{x\}$ when $d_G(x_1)=2$ and $A:= \{x,x_1'\}$ when  $d_G(x_1)=3$.    Let $c: E(G\less A)\rightarrow [5]$ be any star $5$-edge-coloring of $G\less A$.  We may assume that  $ c(uy)=1$ and  $c(uz)=2$. We next prove  that \medskip

\noindent ($*$)   $1 \in c(y_{1})$ and $2\in c(z_1)$.\medskip

  Suppose that  $1 \notin c(y_{1})$ or $2\notin c(z_1)$, say the former.
 If $  c(w) \cup \{1,2\}\ne  [5]$,  then we obtain a star $5$-edge-coloring of $G$  from $c$ by coloring the remaining edges of $G$ as follows (we only consider the worst scenario when  both $x'$ and $x_1'$  exist):
 color the edge $xx_1$ by a color $\alpha$ in $[5] \less (c(w) \cup \{1,2\})$,  $x_1x_1'$ by a color $\beta$ in $[5] \less (c(w) \cup \{\alpha\})$, $ux$ by a color $\gamma$ in $[5] \less  \{1,2, \alpha, c(zz_{1})\}$ and $xx'$ by a color in $[5] \less  \{1,2, \alpha, \gamma\}$, a contradiction. Thus  $ c(w) \cup \{1,2\}= [5]$. Then  $c(w)=\{3,4,5\}$. We may assume that  $c(x_1w)=3$.   If  $c(z) \cup \{1,3\}\ne [5]$, then $\{4,5\}\less c(z)\ne\emptyset$ and we obtain a star $5$-edge-coloring of $G$  from $c$  by coloring  the edge $xx_1$ by  color  $2$, $x_1x_1'$ by color $1$, $ux$ by a color $\alpha$ in $\{4,5\}\less c(z)$   and $xx'$ by a color in  $\{4,5\} \less\alpha$, a contradiction.
 Thus $c(z) \cup \{1,3\} = [5]$ and so $c(z) =\{2,4,5\}$. In particular,  $z'$ must exist. We again obtain a star $5$-edge-coloring of $G$  from $c$  by coloring  $ux,xx', xx_1, x_1x_1'$ by colors $3, c(zz_{1}), 2,1$ in order and then recoloring $uz,zz'$ by colors $c(zz'),2$ in order,
 a contradiction.  Thus $1 \in c(y_{1})$ and $2\in c(z_1)$. This proves ($*$).   \medskip

  By ($*$),   $1 \in c(y_{1})$ and $2\in c(z_1)$.  Then $y_1\ne z_1$,  and $c(yy_1), c(zz_1)\notin \{1,2\}$. We may further assume that    $ c(zz_1)=3$.  Let  $\alpha, \beta\notin  c(z_{1})$ and let   $\gamma, \lambda \notin  c(y_{1})$, where $\alpha, \beta,  \gamma, \lambda\in [5]$.  Since $\alpha, \beta\notin  c(z_{1})$, we may assume that $ c(yy_1)\ne\alpha$.
  We may further assume that $\gamma\ne \alpha$. If $\lambda\ne\alpha$ or $\gamma\notin\{3,\beta\}$,
    then we obtain a star $5$-edge-coloring, say $c'$, of  $G\less A$ from $c$ by recoloring the edges $uz, zz', uy, yy'$ by colors $\alpha, \beta, \gamma, \lambda$, respectively. Then  $c'$ is a star $5$-edge-coloring of  $G\less A$ with $c'(uz)\notin c'(z_1)$,  contrary to ($*$).  Thus  $\lambda=\alpha$ and  $\gamma\in\{3,\beta\}$.  By ($*$), $1\in c(y_1)$ and so  $\alpha=\lambda\ne1$ and  $\gamma\ne1$.
  Let $c'$ be obtained from $c$ by recoloring the edges $uz, zz',  yy'$ by colors $\alpha, \beta, \gamma$, respectively. Then  $c'$ is a star $5$-edge-coloring of  $G\less A$ with $c'(uz)\notin c'(z_1)$, which again contradicts ($*$). \medskip

  This completes the proof of Lemma~\ref{noBad}. \hfill\vrule height3pt width6pt depth2pt\bigskip

\begin{lem}\label{main} For any  $3$-vertex $u$ in $H$ with $N_{H}(u) = \{x, y, z\}$, if  both $x$ and $y$ are bad, then $zx_1, zy_1\notin E(H)$,
and $z$ must be a $3_0$-vertex in $H$, where $x_1$ and $ y_1$ are the other neighbors of $x$ and $y$ in $H$, respectively.
\end{lem}

 \pf   Let $u, x, y, z, x_1, y_1$ be given as in the statement.    Since  $d_{H}(x)=d_{H}(y)=2$, by Lemma~\ref{3nbr}, $u, x, y, z$  are all distinct.   By Lemma~\ref{noBad},  $d_{H}(z)=3$.  Clearly,  both $x_1$ and $y_1$ are bad and so $z\ne x_1, y_1$. By Lemma~\ref{noC3}, $xy\notin E(G)$  and so  $N_{H}(u)$ is an independent set in $H$.
 By Lemma~\ref{noC4},   $x_1\ne y_1$.  It follows that $u, x, y, z, x_1, y_1$ are all distinct.
We first show that $zx_1, zy_1\notin E(H)$.   Suppose that $zx_1\in E(H)$ or  $zy_1\in E(H)$, say the latter.  Then $zy_1$ is not a multiple edge because $d_H(y_1)=2$.    Let   $z_1$ be the third neighbor of $z$ in $H$.  By Lemma~\ref{2nbr}, let  $v$ be the other neighbor of $x_1$ in $H$. Then $v\ne y_1$.  Notice that  $x_1$ and $z_1$  are not necessarily distinct.   Let $A=\{u, x, y, y_1, x_1'\}$. 
Since $G$ is $5$-critical,  let $c: E(G\less A)\rightarrow [5]$ be any star $5$-edge-coloring of $G\less A$.  We may assume that $1,2\notin c(z_1)$ and $c(zz_1)=3$. Let $\alpha \in [5]\less( c(v)\cup \{1\})$ and $\beta\in [5]\less(c(v)\cup\{\alpha\})$. Then we obtain a star $5$-edge-coloring of $G$ from $c$  by  first coloring the  edges $uz, zy_1, xx_1, x_1x_1'$ by colors $1, 2,\alpha, \beta$ in order, and then coloring $ux$ by a color $\gamma$ in $[5]\less \{1,\alpha,\beta,c(x_1v)\}$, $xx'$ by a color in $[5]\less \{1,\alpha,\gamma,c(x_1v)\}$, $uy$ by a color $\theta$ in $[5]\less \{1,2,3, \gamma\}$, $yy_1$ by a color $\mu$ in $[5]\less \{1,2,\gamma,\theta\}$, $yy'$ by a color in $[5]\less \{2,\gamma,\theta,\mu\}$, $y_1y_1'$ by a color in $[5]\less \{1,2, \mu\}$, a contradiction. 
This proves that    $zx_1, zy_1\notin E(H)$.
\medskip

It remains to show that $z$ must be a $3_0$-vertex in $H$.  Suppose that $z$ is not a $3_0$-vertex in $H$. Since $d_{H}(u)=3$, we see that $z$ is either a $3_1$-vertex or a $3_2$-vertex in $H$.
Let $N_{H}(z)=\{u,s,t\}$ with $d_{H}(s)=2$.  By Lemma~\ref{2nbr} applied to the vertex $s$,  $s\ne t$. Since $zx_1, zy_1\notin E(H)$, we see that $x_1, y_1, s,t$ are all distinct. By Lemma~\ref{2nbr}, let $v, w, r$ be the other neighbor of $x_1, y_1, s$ in $H$, respectively.  
Note  that $r$, $t$, $v$, $w$ are not necessarily distinct. By Lemma~\ref{noC4}, both $v$ and $w$ must be  $3$-vertices in $H$. We next prove that  \medskip

\noindent  (a)    if $x'$ or $y'$ exists, then for  any star $5$-edge-coloring $c^*$  of $G\less \{x', y'\}$,  $c^*(xx_1)\in c^*(v)$ or $c^*(yy_1)\in c^*(w)$.\medskip

To see why  (a) is true, suppose that there exists a star  $5$-edge-coloring $c^*: E(G\less \{x', y'\})\rightarrow [5]$ such that $c^*(xx_1)\notin c^*(v)$ and  $c^*(yy_1)\notin c^*(w)$.  Then we obtain a star $5$-edge-coloring of $G$ from $c^*$  by  coloring  $xx'$ by a color in $[5]\less(\{c^*(xx_1)\}\cup c^*(u))$ and $yy'$ by a color in $[5]\less(\{c^*(yy_1)\}\cup c^*(u))$, a contradiction. This proves (a). \medskip

Let $A$ be the set containing $x, y$  and the  $1$-neighbor of  each of $ x_1, y_1$ in $G$ if it exists.    Since $G$ is $5$-critical,  let $c_1: E(G\less A)\rightarrow [5]$ be any star $5$-edge-coloring of $G\less A$.  Let $c$ be a star $5$-edge-coloring of $G\less\{x,  x', y', x_1'\}$ obtained from $c_1$ by  coloring  $yy_1$ by a color $\alpha$ in $[5]\less (c_1(w)\cup\{c_1(uz)\})$,    $uy$ by a color  in $[5]\less(c_1(z)\cup \{\alpha\})$, and  $y_1y_1'$  by a color $\beta$ in $[5]\less (c_1(w)\cup \{\alpha\})$.  We may assume that $c(uz)=1$, $c(zs)=2$ and $c(zt)=3$. By the choice of $c(uy)$, we may further assume that
$c(uy)=4$. We next obtain a contradiction  by extending $c$ to be a   star  $5$-edge-coloring of $G$  (when neither of  $x'$ and $y'$  exists) or a star $5$-edge-coloring of $G\less\{x', y'\}$ (when $x'$ or $y'$ exists)  which violates (a).  We consider the worst scenario when $x'$ and  $y'$ exist.  We first prove two claims. \bigskip

\noindent  {\bf Claim 1}:  $\beta=4$ or  $c(y_1w)=4$.\medskip


 \pf Suppose that $\beta\ne 4$ and $c(y_1w)\ne4$.  We next show  that $c(v)\cup\{1,4\}\ne [5]$. Suppose that $c(v)\cup\{1,4\}= [5]$. Then $c(v)=\{2,3,5\}$. Clearly,   $c(x_1v)=5$, otherwise, coloring $ux$, $xx_1$,   $x_1x'_1$ by colors  $5,1,4$ in order,  we obtain    a star $5$-edge-coloring of $G\less\{x', y'\}$  which violates (a), a contradiction.   
We see that   $1\in c(s)\cap c(t)$, otherwise,  we may assume  that   $1\notin c(s)$, we obtain   a star $5$-edge-coloring of $G\less\{x', y'\}$  which violates (a) as follows:  when  $\alpha\ne2$, color $ux, xx_1, x_1x_1'$  by colors $2, 4,1$ in order; when $\alpha=2$,  first color $ux, xx_1, x_1x_1'$ by  colors $2, 4,1$ in order and then  recolor  $yy_1, y_1y_1'$ by colors $\beta,2$ in order.  It follows that  $4, 5\in c(s)\cup c(t)$, otherwise,   say $\theta\in\{4,5\}$ is not in $ c(s)\cup c(t)$,  let $\alpha'\in\{2,3\}\less \alpha$, now either  coloring $ux, xx_1, x_1x_1'$ by colors $\alpha', 4,1$ in order and  then recoloring  $uz$ by color $5$ when $\theta=5$; or  coloring $ux, xx_1, x_1x_1'$ by colors $\alpha', 1, 4$ in order and then recoloring  $uz, uy$ by colors $4, 1$ in order when $\theta=4$,
 we obtain  a star $5$-edge-coloring of $G\less\{x', y'\}$  which violates (a).  Thus $c(s)=\{1,2,\theta\}$ and   $c(t)=\{1,3,9-\theta\}$, where $\theta\in\{4,5\}$.
  If  $c(ss')=\theta$  or $c(sr)= \theta$  and $2\notin c(r)$, then we obtain  a star $5$-edge-coloring of $G\less\{x', y'\}$ (which violates (a)) as follows: when $\theta=5$,  color $ux, xx_1, x_1x_1'$ by colors $3,1,4$ in order  and then recolor  $uz$ by color $5$; when $\theta=4$ and $\alpha\in\{2,5\}$, first  color  $ux, xx_1, x_1x'_1$   by colors $3, 1, 4$ in order, and then   recolor  $uz, uy$ by colors $4, 1$ in order;
  when $\theta=4$ and $\alpha=3$ and $\beta\ne 5$,  color  $ux, xx_1, x_1x'_1$   by colors $5, 1, 4$ in order and then   recolor  $uz, uy, yy_1, y_1y_1'$ by colors $4, 3, \beta, 3$ in order;  when $\theta=4$ and $\alpha=3$ and $\beta=5$,  color  $ux, xx_1, x_1x'_1$   by colors $3, 1, 4$ in order and then   recolor  $uz, uy, yy_1, y_1y_1'$ by colors $4, 1, 5, 3$ in order.
 Thus   $c(ss')=1$, $c(sr)=\theta$ and $2\in c(r)$.  Now recoloring the edge $ss'$ by a color in $\{3,9-\theta\}\less c(r)$  yields a star $5$-edge-coloring $c$ of $G\less\{x,  x', y', x_1'\}$ satisfying $\beta\ne 4$, $c(y_1w)\ne4$,  $c(v)\cup\{1,4\}= [5]$ and $c(x_1v)=5$ but   $1\notin c(s)\cap c(t)$, a contradiction. This proves that  $c(v)\cup\{1,4\}\ne [5]$. \medskip

Since $c(v)\cup\{1,4\}\ne [5]$, we see that $[5]\less (c(v)\cup\{1,4\})= \{5\}$, otherwise, coloring $ux$ by color $5$,  $xx_1$ by a color $\gamma$ in $[5]\less(c(v)\cup\{1,4,5\})$, and  $x_1x'_1$ by  a color in $[5]\less (c(v)\cup \gamma)$, we obtain a star $5$-edge-coloring of $G\less\{x', y'\}$ which violates (a).  Clearly, $2,3\in c(v)$ and $\{1,4\}\less c(v)\ne\emptyset$. Let $\gamma\in \{1,4\}\less c(v)$ and   $\alpha'\in\{2,3\}\less \alpha$.
Then   $1\in c(s)\cap c(t)$, otherwise,  we may assume  that   $1\notin c(s)$, now      coloring $ux, xx_1, x_1x_1'$ by colors $2, 5,\gamma$ in order  yields   a star $5$-edge-coloring of $G\less\{x', y'\}$  which violates (a).
It follows that  $4, 5\in c(s)\cup c(t)$, otherwise, say $\theta\in\{4,5\}$ is not in $ c(s)\cup c(t)$,  first recoloring $uz$ by color $\theta$ and then  either  coloring $ux, xx_1, x_1x_1'$ by colors $\alpha', 5,  \gamma$ in order and then recoloring $uy$ by color $1$ when $\theta=4$; or coloring
$ux, xx_1, x_1x_1'$ by colors $\alpha', 1, 5$ in order when $\theta=5$ and $\gamma=1$; or
coloring $ux, xx_1, x_1x_1'$ by colors $1,4,5$ in order  when  $\theta=5$, $\gamma=4$ and $c(x_1v)\ne1$; or coloring $ux, xx_1, x_1x_1'$ by colors $\alpha',4,5$ in order  when  $\theta=5$, $\gamma=4$ and $c(x_1v)=1$,
we obtain a star $5$-edge-coloring of $G\less\{x', y'\}$  which violates (a).
Thus $c(s)=\{1,2,\theta\}$ and   $c(t)=\{1,3,9-\theta\}$, where $\theta\in\{4,5\}$. If  $c(ss')=\theta$  or $c(sr)= \theta$  and $2\notin c(r)$, then we obtain  a star $5$-edge-coloring of $G\less\{x', y'\}$ (which violates (a)) as follows: when $\theta=5$ and $\gamma=1$,  color $ux, xx_1, x_1x_1'$ by colors $3, 1, 5$ in order and then recolor  $uz$ by colors $5$;
when $\theta=5$,  $\gamma=4$ and $c(x_1v)\ne1$,  color $ux, xx_1, x_1x_1'$ by color $ 1, 4, 5$ in order and then recolor  $uz$ by colors $5$;
when $\theta=5$,  $\gamma=4$ and $c(x_1v)=1$,  color $ux, xx_1, x_1x_1'$ by color $ 3, 4, 5$ in order and then recolor  $uz$ by colors $5$ (and  further recolor $yy_1$ by $\beta$ and $y_1y_1'$ by $\alpha$ when $\alpha=3$);
when $\theta=4$ and $\beta\ne1$,  color $ux, xx_1, x_1x_1'$ by color $3, 5, \gamma$ in order and then recolor  $uz,  uy$ by colors  $4, 1$ in order,   and finally recolor $yy_1$ by a color  $\beta'\in \{\alpha,\beta\}\less 3$ and  $y_1y'_1$ by a color in $\{\alpha,\beta\}\less \beta'$;
when $\theta=4$,  $\beta=1$ and  $\gamma=1$,  color $ux, xx_1, x_1x_1'$ by color $5,1,5$ in order and then recolor  $uz,  uy, yy_1, y_1y_1'$ by colors  $4, 3, 1, \alpha$ in order;
when $\theta=4$,  $\beta=1$,   $\gamma=4$ and $\alpha\ne3$,  color $ux, xx_1, x_1x_1'$ by color $3, 5, 4$ in order and then recolor  $uz,  uy$ by colors  $4, 1$ in order;
when $\theta=4$,  $\beta=1$,   $\gamma=4$ and $\alpha=3$, let $\gamma'\in\{1,3\}\less c(x_1v)$,  color $ux, xx_1, x_1x_1'$ by color $\gamma', 5, 4$ in order and then recolor  $uz$ by color $4$,  $uy$ by color  $5$,  $yy_1$ by a color $\beta'$ in $\{1,3\}\less \gamma'$ and  $y_1y_1'$ by a color in $\{1,3\}\less \beta'$.
Thus $c(ss')=1$, $c(sr)=\theta$ and $2\in c(r)$.  Now recoloring the edge $ss'$ by a color in $\{3,9-\theta\}\less c(r)$  yields a star $5$-edge-coloring $c$ of $G\less\{x,  x', y', x_1'\}$    satisfying $\beta\ne 4$,  $c(y_1w)\ne4$  and $[5]\less (c(v)\cup\{1,4\})= \{5\}$ but   $1\notin c(s)\cap c(t)$, a contradiction.  This completes the proof of Claim 1.
 \hfill\vrule height3pt width6pt depth2pt\\

 \noindent  {\bf Claim 2}:  $\beta=4$.  \medskip

Suppose that  $\beta\ne 4$. By Claim 1,  $c(y_1w)=4$.  We first consider the case when  $c(w)=\{2,3,4\}$. Then  $\alpha=5$ and $\beta=1$. We claim that
$c(v)\cup\{1,4\}\ne [5]$. Suppose that $c(v)\cup\{1,4\}= [5]$.
Then $c(v)=\{2,3,5\}$. Clearly,  $1\in c(s)\cap c(t)$, otherwise,  we may assume that $1\notin c(s)$, now coloring $ux, x x_1, x_1x'_1$ by colors $5,4,1$ in order and then  recoloring $uy$ by $2$,  we  obtain a star $5$-edge-coloring of $G\less\{x', y'\}$ which violates (a).
 It follows that  $4, 5\in c(s)\cup c(t)$, otherwise, say $\theta\in\{4,5\}$ is not in $ c(s)\cup c(t)$,  now coloring $ux, x x_1, x_1x'_1$ by colors $3,1,4$ in order and then  recoloring $uz,  uy, yy_1, y_1y_1'$ by colors $\theta, 2, 1,5$ in order
we obtain a star $5$-edge-coloring of $G\less\{x', y'\}$  which violates (a).
Thus $c(s)=\{1,2,\theta\}$ and   $c(t)=\{1,3,9-\theta\}$, where $\theta\in\{4,5\}$.
If  $c(ss')=\theta$  or $c(sr)= \theta$  and $2\notin c(r)$, then   coloring $ux, x x_1, x_1x'_1$ by colors $3,1,4$ in order and then  recoloring $uz,  uy, yy_1, y_1y_1'$ by colors $\theta, 9-\theta, 1,5$ in order yileds a star $5$-edge-coloring of $G\less\{x', y'\}$ which violates (a).
Thus $c(ss')=1$,  $c(sr)= \theta$  and $2\in c(r)$. Now recoloring the edge $ss'$ by a color in $\{3,9-\theta\}\less c(r)$  yields a star $5$-edge-coloring $c$ of $G\less\{x,  x', y', x_1'\}$  satisfying $\alpha=5$, $\beta=1$, $c(y_1w)=4$ and  $c(v)\cup\{1,4\}= [5]$ but $1\notin c(s)\cap c(t)$, a contradiction. This proves that  $c(v)\cup\{1,4\}\ne [5]$.
 Let $\eta=5$ when $5\notin c(v)$ or $\eta \in \{2,3\}\less c(v)$   when $5\in c(v)$.  Let $\mu\in [5]\less( c(v)\cup\{\eta\})$.  By Claim 1 and the symmetry between $x$ and $y$, either $4\notin c(v)$ or $5\notin c(v)$. We see that   $\mu=4$ when $\eta\ne 5$.
  Then  $1\in c(s)\cap c(t)$, otherwise, we may assume $1\notin c(s)$, we obtain a star $5$-edge-coloring of $G\less\{x', y'\}$  (which violates (a)) as follows: when  $\eta\neq 2$,  color $ux, xx_1, x_1x_1'$ by colors $2, \eta, \mu$ in order; when $\eta=2$, then $\mu=4$,  first recolor $uy$ by color $2$ and  then color $ux, xx_1, x_1x'_1$ by colors  $5,4,  2$ in order.
It follows that  $4, 5\in c(s)\cup c(t)$, otherwise, say $\theta\in\{4,5\}$ is not in $ c(s)\cup c(t)$,
now first recoloring $uz, yy_1, y_1y_1'$ by colors $\theta, 1,5$  in order, and then coloring  $xx_1, x_1x'_1$ by colors $\eta,\mu$ in order,  $ux$ by a color $\gamma$ in $[5]\less\{\mu,\eta,\theta, c(x_1v)\}$,  and finally  coloring $uy$ either by a color in $\{2,3\}\less\eta$ when  $\gamma=1$ or by a color in  $\{2,3\}\less\gamma$ when $\gamma\ne1$,
we obtain a star $5$-edge-coloring of $G\less\{x', y'\}$  which violates (a).
Thus $c(s)=\{1,2,\theta\}$ and   $c(t)=\{1,3,9-\theta\}$, where $\theta\in\{4,5\}$.
If   $c(ss')=\theta$  or $c(sr)= \theta$  and $2\notin c(r)$, we obtain a star $5$-edge-coloring of $G\less\{x', y'\}$  (which violates (a)) as follows:
when $\theta=4$ and $\eta=5$, color $ux, xx_1, x_1x_1'$ by colors $3, 5, \mu$ in order and then recolor $uz, uy$ by colors $4,1$ in order;
when $\theta=4$ and $\eta\in\{2,3\}$, then $\mu=4$,  first recolor $uz, uy$ by colors $4,3$ in order and then color $xx_1, x_1x_1'$ by colors $ \eta, 4$ in order and finally color $ux$ by a color $\gamma$ in $\{1,5\}\less c(x_1v)$, $yy_1$ by a color $\lambda$ in $\{1,5\}\less \gamma$, and $y_1y_1'$ by a color in $\{1,5\}\less \lambda$;
when $\theta=5$ and $\eta\in\{2,3\}$, then $\mu=4$, color $ux, xx_1, x_1x_1'$ by colors $1, 4, \eta$ in order and then recolor $uz, uy, yy_1, y_1y_1'$ by colors $5,3, 1,5$ in order;
when $\theta=5$,  $\eta=5$ and $\mu\ne3$, color $ux, xx_1, x_1x_1'$ by colors $1, \mu, 5$ in order and then recolor $uz, uy, yy_1, y_1y_1'$ by colors $5,3, 1,5$ in order;
when $\theta=5$,  $\eta=5$ and $\mu=3$, first recolor $uz, uy, yy_1, y_1y_1'$ by colors $5,3, 1,5$ in order, then color $xx_1, x_1x_1'$ by colors $5, 3$ in order and  finally color $ux$ by a color in $\{1,4\}\less c(x_1v)$.
 Thus $c(ss')=1$, $c(sr)= \theta$  and $2\in c(r)$.  Now recoloring the edge $ss'$ by a color in $\{3,9-\theta\}\less c(r)$  yields a star $5$-edge-coloring $c$ of $G\less\{x,  x', y', x_1'\}$  satisfying $\alpha=5$, $\beta=1$, $c(z)=\{1,2,3\}$, $c(uy)=c(y_1w)=4$ and  $c(v)\cup\{1,4\}\ne[5]$ but $1\notin c(s)\cap c(t)$, a contradiction.\\

We next consider the case when  $c(w)\neq\{2,3,4\}$.
If $\alpha, \beta\ne5$, then recoloring $uy$ by color $5$ yields  a  star $5$-edge-coloring $c$ of $G\less\{x,  x', y', x_1'\}$ with $c(uy)\ne c(y_1y_1'), c(y_1w)$,  contrary to  Claim 1. Thus either $\alpha=5$ or  $\beta=5$.  Then $1\in c(w)$ because $c(w)\neq\{2,3,4\}$ and $|c(w)|=3$.
It follows that $\alpha, \beta\in\{2,3,5\}$ and $5\in\{\alpha, \beta\}$. We may assume that $\alpha\in\{2,3\}$ and $\beta=5$ by permuting the colors on $yy_1$ and $y_1y_1'$ if needed.  Then $4,5\in c(s)\cup c(t)$, otherwise, say $\theta\in\{4,5\}$ is not in $ c(s)\cup c(t)$, we obtain a a star $5$-edge-coloring $c$ of $G\less\{x,  x', y', x_1'\}$  which contradicts Claim 1 by recoloring $uz, uy$ by colors $\theta, 1$ in order.  Let $\alpha'\in\{2,3\}\less \alpha$.  We next show that $c(ss')=1$,  $c(sr)= \theta$  and $2\in c(r)$.\medskip

Suppose first  that $c(v)\cup\{1,4\}= [5]$.
Then $c(v)=\{2,3,5\}$. We see that  $c(x_1v)=5$, otherwise, coloring $ux$, $xx_1$, $x_1x'_1$ by colors  $5,1,4$ in order,  we obtain    a star $5$-edge-coloring of $G\less\{x', y'\}$  which violates (a).  Clearly,  $1\in c(s)\cap c(t)$, otherwise,  we may assume that $1\notin c(s)$, now coloring $ux, x x_1, x_1x'_1$ by colors $2,4,1$ in order and then  recoloring $yy_1,y_1y_1'$ by colors $5,\alpha$, we  obtain a star $5$-edge-coloring of $G\less\{x', y'\}$ which violates (a).  Since $4,5\in c(s)\cup c(t)$, we see that  $c(s)=\{1,2,\theta\}$ and   $c(t)=\{1,3,9-\theta\}$, where $\theta\in\{4,5\}$.
If  $c(ss')=\theta$  or $c(sr)= \theta$  and $2\notin c(r)$, then  recoloring $uz, uy$ by colors $\theta, 1$ in order  yields  a  star $5$-edge-coloring $c$ of $G\less\{x,  x', y', x_1'\}$ with $c(uy)\ne c(y_1y_1'), c(y_1w)$,  contrary to  Claim 1.
Thus $c(ss')=1$,  $c(sr)= \theta$  and $2\in c(r)$.  Next suppose   that $c(v)\cup\{1,4\}\ne  [5]$.
Let $\eta=5$ when $5\notin c(v)$ or $\eta \in \{2,3\}\less c(v)$   when $5\in c(v)$.  Let $\mu\in [5]\less( c(v)\cup\{\eta\})$. By Claim 1 and the symmetry between $x$ and $y$, either $4\notin c(v)$ or $5\notin c(v)$. We see that   $\mu=4$ when $\eta\ne 5$.   Then  $1\in c(s)\cap c(t)$, otherwise, we may assume $1\notin c(s)$, we obtain a star $5$-edge-coloring of $G\less\{x', y'\}$  (which violates (a)) as follows: when  $\eta=5$,  color $ux, xx_1, x_1x_1'$ by colors $4, 5, \mu$ in order and then recolor $uy, yy_1, y_1y_1'$ by colors $2,5,\alpha$ in order; when $\eta\in\{2,3\}$, then $\mu=4$,  color $ux, xx_1, x_1x'_1$ by colors  $5,\eta,  4$ in order.  Since $4,5\in c(s)\cup c(t)$, we see that  $c(s)=\{1,2,\theta\}$ and   $c(t)=\{1,3,9-\theta\}$, where $\theta\in\{4,5\}$.
If  $c(ss')=\theta$  or $c(sr)= \theta$  and $2\notin c(r)$, then  recoloring $uz, uy$ by colors $\theta, 1$ in order  yields  a  star $5$-edge-coloring $c$ of $G\less\{x,  x', y', x_1'\}$ with $c(uy)\ne c(y_1y_1'), c(y_1w)$,  contrary to  Claim 1.
Thus $c(ss')=1$,  $c(sr)= \theta$  and $2\in c(r)$. \medskip

Now recoloring the edge $ss'$ by a color in $\{3,9-\theta\}\less c(r)$  yields a star $5$-edge-coloring $c$ of $G\less\{x,  x', y', x_1'\}$  satisfying $\alpha\in\{2,3\}$, $\beta=5$, $c(y_1w)=4$ and   $c(w)\neq\{2,3,4\}$ but $1\notin c(s)\cap c(t)$, a contradiction.   This completes the proof of Claim 2. \hfill\vrule height3pt width6pt depth2pt\\


By Claim 2, $\beta=4$.     Suppose that $\alpha\ne5$.  Then  $\alpha\in\{2,3\}$.  Note that $\alpha\notin c(w)\cup\{1\}$. Now  recoloring $uy$ by color $5$, we obtain a star $5$-edge-coloring $c$ of $G\less\{x,  x', y', x_1'\}$ satisfying $c(uz)=1$, $c(zs)=2$ and   $c(zt)=3$ but $\beta\ne c(uy)$, contrary to Claim 2. Thus $\alpha=5$ and   so $ c(w)=\{1,2,3\}$.
 By the symmetry of $x$ and $y$, $c(v)=\{1,2,3\}$.  Then $1\in c(s)\cap c(t)$, otherwise,   we may assume  that   $1\notin c(s)$, now      coloring $ux, xx_1, x_1x_1'$ by colors $2,5,4$ in order  yields   a star $5$-edge-coloring of $G\less\{x', y'\}$  which violates (a).  It follows that  $4, 5\in c(s)\cup c(t)$, otherwise,   say $\theta\in\{4,5\}$ is not in $ c(s)\cup c(t)$,   now first coloring $ux, xx_1, x_1x_1'$ by colors $2,9-\theta,  \theta$ in order and then recoloring $uz, uy, yy_1, y_1y_1'$ by colors $\theta, 3, 9-\theta, \theta$ in order, we obtain      a star $5$-edge-coloring of $G\less\{x', y'\}$  which violates (a).  Thus $c(s)=\{1,2,\theta\}$ and   $c(t)=\{1,3,9-\theta\}$, where $\theta\in\{4,5\}$.  If $c(ss')=\theta$ or $c(sr)=\theta$ and $2\notin c(r)$, then we obtain  a star $5$-edge-coloring of $G\less\{x', y'\}$ (which violates (a)) by coloring $ux, xx_1, x_1x_1'$ by colors $1, 9-\theta, \theta$ in order, and then  recoloring $uz, uy, yy_1, y_1y_1'$ by colors $\theta, 3, 9-\theta, \theta$ in order.
Thus $c(ss')=1$ and   $2\in c(r)$. Now recoloring $ss'$ by a color in $\{3,9-\theta\}\less c(r)$,  we obtain a star $5$-edge-coloring $c$ of $G\less\{x,  x', y', x_1'\}$ satisfying $c(uz)=1$, $c(zs)=2$,  $c(zt)=3$, $\beta= 4$ and  $\alpha=5$  but   $1\notin c(s)\cap c(t)$. \medskip

This completes the proof of Lemma~\ref{main}.\hfill\vrule height3pt width6pt depth2pt\\


 \section{Proof of Theorem~\ref{main}}
We are now ready to prove Theorem~\ref{main}.  Suppose the  assertion is false. Let $G$ be a subcubic multigraph with $\text{mad}(G)<12/5$ and  $\chi'_s(G)>5$. Among all counterexamples we choose $G$ so that
$|G|$ is minimum. By the choice of $G$, $G$ is connected,  star $5$-critical,  and  $\text{mad}(G)<12/5$.  For all $i\in[3]$, let $A_i=\{v\in V(G): \, d_G(v)=i\}$ and let  $n_i=|A_i|$ for all $i\in[3]$.
Since  $\text{mad}(G)<12/5$, we see that  $3n_3<2n_2+7n_1$ and so  $A_1\cup A_2\ne\emptyset$. By Lemma~\ref{deg=1}(a), $A_1$  is an  independent set in $G$ and $N_G(A_1)\subseteq A_3$.  Let $H=G\less A_1$.  Then $H$ is connected and $\text{mad}(H)<12/5$. By Lemma~\ref{deg=1}(b), $\delta(H)\ge2$.  By Lemma~\ref{3nbr}, every  $3_2$-vertex   in $H$ has three distinct neighbors in $H$.
We say that a $3_2$-vertex   in $H$ is \dfn{bad} if  both of its  $2$-neighbors  are bad.  A vertex $u$ is a \dfn{good} (resp. \dfn{bad}) $2$-neighbor of a vertex $v$ in $H$  if $uv\in E(H)$ and $u$ is a good (resp. bad) $2$-vertex. By Lemma~\ref{main},  every bad $3_2$-vertex in $H$ has a unique $3_0$-neighbor. We now  apply the discharging method to obtain a contradiction. \medskip

For each vertex $v\in V(H)$, let $\omega(v):= d_{H}(v)-\frac{12}{5}$ be the initial charge of $v$. Then $ \sum_{v\in V(H)} \omega(v) =2e(H)-\frac{12}{5}|H|=|H|(2e(H)/|H|-\frac{12}{5})<0$. Notice that for each $v\in V(H)$,  $\omega(v)=2-\frac{12}{5}=-\frac{2}{5}$ if $d_{H}(v)=2$,  and $\omega(v)=3-\frac{12}{5}=\frac{3}{5}$ if $d_{H}(v)=3$. We will redistribute the  charges of vertices in  $H$ as follows. \medskip

\myitemitem  {(R1):}  every bad $3_2$-vertex in $H$ takes  $\frac1 {5}$ from its unique $3_0$-neighbor.

\myitemitem {(R2):} every $3_1$-vertex in $H$ gives $\frac3 {5}$ to its  unique $2$-neighbor.

\myitemitem {(R3):}  every  $3_2$-vertex  in $H$ gives  $\frac1 {5}$ to each of its good $2$-neighbors (possibly none) and $\frac25$ to each of its bad $2$-neighbors (possibly none).

\myitemitem {(R4):} every $3_3$-vertex in $H$ gives $\frac1 {5}$ to each of its $2$-neighbors.\medskip

 Let $\omega^*$ be the new charge of $H$ after applying the above discharging rules in order. It suffices to show that $\sum_{v\in V(H)} \omega^*(v)\geq0$.
For any $v\in V(H)$ with $d_{H}(v)=2$, by Lemma~\ref{2nbr}, $v$ has two distinct neighbors in $H$.
If $v$ is a good $2$-vertex, then $v$ takes  at least $\frac15$ from each of its $3$-neighbors under  (R2), (R3) and (R4), and so  $\omega^*(v)\geq0$.  Next, if  $v$ is a bad $2$-vertex, let $x$, $y$ be the  two neighbors of $v$ in  $H$. We may assume that  $y$ is a bad $2$-vertex. By Lemma~\ref{2nbr}, let $z$ be the other  neighbor  of $y$ in $H$.  By Lemma~\ref{noC4}, we may assume that $d_{H}(x)=3$.  By Lemma~\ref{noBad}, $x$ is either a $3_1$-vertex or a $3_2$-vertex in $H$.  Under (R2) and (R3), $v$ takes at least $\frac25$ from $x$.  If $d_{H}(z)=3$, then  by a similar argument,  $y$ must take at least $\frac25$ from $z$.  In this case, $\omega^*(v)+\omega^*(y)\geq0$.  If $d_{H}(z)=2$, then $z$ is bad. By Lemma~\ref{2nbr}, let $w$ be the other neighbor of $z$.  By Lemma~\ref{noC4}, each of $x$ and $w$ must be a $3_1$-vertex in $H$. Under (R2), $v$ takes $\frac35$ from $x$ and $z$ takes $\frac35$ from $w$. Hence, $\omega^*(v)+\omega^*(y)+\omega^*(z)\geq0$. \medskip

For any $v\in V(H)$ with $d_{H}(v)=3$, if $v$ is a bad $3_2$-vertex,  then   $v$ has a unique $3_0$-neighbor by Lemma~\ref{main}. Under (R1) and (R3), $v$ first takes $\frac15$ from its unique $3_0$-neighbor and then gives $\frac25$ to each of its bad $2$-neighbors, we see that  $\omega^*(v)\geq0$.  If  $v$ is not a bad $3_2$-vertex, then  $v$ gives either nothing  or one of $\frac15$,  $\frac25$,  and $\frac35$ in total to its neighbors under    (R1), (R2), (R3) and (R4). In either case,   $\omega^*(v)\ge0$.
Consequently,  $\sum_{v\in V(H)} \omega^*(v)\geq0$, contrary to the fact that  $\sum_{v\in V(H)} \omega^*(v)=\sum_{v\in V(H)}\omega(v)<0$. \medskip

This completes the proof of Theorem~\ref{main}.
\hfill\vrule height3pt width6pt depth2pt\medskip

  \vspace{5mm} \noindent {\bf Acknowledgments.}  Zi-Xia Song would like to thank Yongtang Shi and  the Chern Institute of Mathematics at Nankai University for hospitality and support during her visit  in May 2017. \medskip 
  
 \noindent     Hui  Lei and Yongtang  Shi are partially supported by the National Natural Science Foundation of China and the Natural Science Foundation of Tianjin (No.17JCQNJC00300). \medskip
  
 \noindent  Tao  Wang is partially supported by the National Natural Science Foundation of China (11101125) and the Fundamental Research Funds for Universities in Henan (YQPY20140051).

\end{document}